\documentclass[12pt]{amsart}
 \usepackage{graphics}

\begin{document}
 
\newtheorem{lemma}{Lemma}[section]
\newtheorem{prop}[lemma]{Proposition}
\newtheorem{cor}[lemma]{Corollary}
\newtheorem{thm}{Theorem}

\theoremstyle{definition}
\newtheorem{rem}[lemma]{Remark}
\newtheorem{rems}[lemma]{Remarks}
\newtheorem{defi}[lemma]{Definition}
\newtheorem{ex}{Example}
\newtheorem{convention}[lemma]{Convention}

\newcommand{\cH}{\mathcal H}
\newcommand{\oh}{\overline{h}}
\newcommand{\ok}{\overline{k}}
\title[Walls and cubings]{From wall spaces to CAT(0) cube complexes}
\author{Indira Chatterji$^\dagger$ and Graham Niblo}
\email{indira@math.cornell.edu, G.A.Niblo@maths.soton.ac.uk}
\date{\today}
\thanks{$^\dagger$ Partially supported by the Swiss National Funds}
\begin{abstract}We explain how to adapt a construction due to M.~Sageev in order to construct a proper action of a group on a CAT(0) cube complex starting from a proper action of the group on a wall space.\end{abstract}
\maketitle
\section*{Introduction}
A group is said to be \emph{a-T-menable} if it admits a proper isometric action on a Hilbert space. The existence of such an action yields information about the unitary representation theory of the group, and in particular the Baum-Connes conjecture is known to hold for a-T-menable groups, see \cite{HK}.

In \cite{NR} the second author, together with Martin Roller, described a method to construct isometric actions on Hilbert space by adapting a technique of Sageev \cite{S}. The construction can be interpreted as giving a recipe for embedding a CAT(0) cube complex in the unit cube in a Hilbert space so that the group of cellular isometries of the cube complex acts by affine isometries on the Hilbert space. A similar construction can be used to derive affine actions on Hilbert space from actions on so-called spaces with walls (see \cite{HP}, \cite{welches} or \cite{NR}). In \cite{V} the notion of a space with walls is generalized to that of a measured wall space. The relationship between spaces with walls and measured wall spaces is roughly analogous to the relationship between a classical tree and an ${\bf R}$-tree, and in \cite{V}, Cherix, Martin and Valette showed that a finitely generated group admits a proper action on a measured wall space if and only if it is a-T-menable. It is natural to ask whether or not groups which act properly on measured wall spaces must also admit proper actions on spaces with walls. The motivation of this note is to highlight a partial answer to the following question:
\begin{center}\emph{``Are there discrete a-T-menable groups which cannot\\ act properly on a discrete wall space?''}\end{center} 
We will show that a group acts (properly) isometrically on a space with walls if and only it acts (properly) isometrically on a CAT(0) cube complex of dimension $k$, where $k$ is the cardinality of a maximal family of pairwise crossing walls. The first author with Kim Ruane in \cite{CR} prove that a group acting properly and with uniformly bounded point stabilisers on a CAT(0) cube complex of finite dimension cannot contain amenable subgroups of super-polynomial growth, so we obtain the following.
\begin{thm}\label{principal}Let $G$ be a group containing an amenable subgroup $H$ of super-polynomial growth. If $G$ acts properly on a space with walls and with uniformly bounded point stabilisers then there are arbitrarily large families of walls in the space which pairwise cross.\end{thm}

It should be noted that the hypothesis that stabilisers are uniformly bounded is necessary. In \cite{CR} an  example is given of a group containing an amenable subgroup $H$ of super-polynomial growth which acts properly on a space with walls in which at most two walls pairwise cross. The example has been pointed out by S. Mozes, and consists of $PGL_{2}(F_{p}[t,t^{-1}])$, which acts properly on a product of two trees.

\medskip

We would like to mention that an attractive alternative of the cubulation procedure has been simultaneously and independently described by Bogdan Nica in \cite{Nicas}.
\section{Basic definitions}\label{basics}
In this section we shall recall some basic definitions and notations concerning wall spaces and CAT(0) cube complexes.
A \emph{space with walls} (or \emph{wall space}) is a set $Y$ with a non-empty collection of non-empty subsets $\cH\subset{\mathcal P}(Y)$ called \emph{half-spaces}, closed under the involution
\begin{eqnarray*}*:\cH & \to & \cH\\
h & \mapsto & h^c=Y\setminus h\end{eqnarray*}
and satisfying the following rule (strong finite interval condition):

\medskip

For any $p,q\in Y$, there exists only finitely many half-spaces $h\in\cH$ such that $p\in h$ and $q\in h^c$ (such a half-space $h$ is then said to \emph{separate} $p$ from $q$).

\medskip

The \emph{walls} are the pairs $\oh=\{h,h^c\}$, so that the set of walls is $W=\cH/*$. We say that two walls $\oh$ and $\ok$ \emph{cross} if all four intersections
$$h\cap k,\  h\cap k^c,\ h^c\cap k\hbox{ and }h^c\cap k^c$$
are nonempty. We define the \emph{intersection number} of $W$ by
$$I(W)=\sup\{n\in{\bf N}\cup\infty|\exists \oh_1,\dots,\oh_n\in W\hbox{ pairwise crossing}\}.$$
For two walls $\oh=\{h,h^c\}$ and $\ok=\{k,k^c\}$ that don't cross and $p\in Y$ such that $p\in h\cap k$, we say that $\ok$ \emph{separates $p$ from the wall} $\oh$ if $k\subset h$, namely if the following situation occurs:
\begin{center}
\setlength{\unitlength}{1cm}
\begin{picture}(5,3)
\put(1,0){\line(0,1){2}}
\put(1,2.5){\makebox(0,0){$\oh$}}
\put(1,1){\vector(1,0){1}}
\put(1.3,1.3){\makebox(0,0){$h$}}
\put(3,0){\line(0,1){2}}
\put(3,2.5){\makebox(0,0){$\ok$}}
\put(3,1){\vector(1,0){1}}
\put(3.3,1.3){\makebox(0,0){$k$}}
\put(5,1){\circle*{0.1}}\
\put(5,1.3){\makebox(0,0){$p$}}
\end{picture}\end{center}

Let $G$ be a group acting on the space $Y$, we say that $G$ \emph{acts on the space with walls} if the action is such that
$$g(h)\in\cH\hbox{ for all }h\in\cH\hbox{ and }g\in G.$$
The space with walls $Y$ inherits a pseudo-metric by
$$d(p,q)=\sharp\{\oh=\{h,h^c\}\in W |h\hbox{ separates }p\hbox{ from }q\}.$$
We say that the action is \emph{proper} if it is metrically proper for the above given pseudo-metric (i.e. for any sequence $g_n$ of elements in $G$ tending to infinity, then $d(g_n(p),p)$ tends to infinity in ${\bf R}_+$ for any $p\in Y$). If $G$ acts on the space with walls $Y$, then $G$ also acts on the set of walls $W$. 

Let us now briefly recall that a \emph{cube complex} $X$ is a metric polyhedral complex in which each cell is isomorphic to the Euclidean cube $[0,1]^n$ and the gluing maps are isometries. A cube complex is called CAT(0) if the metric induced by the Euclidean metric on the cubes turns it into a CAT(0) metric space. We shall denote by $X^i$ the $i$-skeleton  of $X$, and say that $X$ is finite dimensional if there is $n<\infty$ such that $X^m$ is empty for any $m>n$. In the sequel we shall use an equivalent characterization of CAT(0) cube complexes, given by the following.
\begin{thm}[Gromov, see Thm 5.4 page 206, \cite{bible}]\label{Gromov}Let $X$ be a cubical complex. Then $X$ is CAT(0) if and only if
\begin{itemize}
\item[(i)]The link of each vertex is a flag complex (i.e. $X$ is locally CAT(0)).
\item[(ii)]$X$ is simply connected.
\end{itemize}\end{thm}
Recall that a \emph{flag complex} is a simplicial complex in which each subgraph isometric to the 1-skeleton of a $k$-dimensional simplex actually is the 1-skeleton of a $k$-dimensional simplex.
\section{Sageev's construction}\label{sageev}
In his thesis Sageev showed how to construct a CAT(0) cube complex starting with a finitely generated group $G$ and a so-called codimension-one subgroup. Here we will adapt Sageev's construction to the more general situation of a group acting on a space with walls. More precisely we explain the proof of the following.
\begin{thm}Let $G$ be  a discrete group acting properly on a space with walls $Y$ and let $W$ denote the set of walls. Then there exists a CAT(0) cube complex $X$ on which $G$ acts properly. Moreover, in the case where $I(W)$ is finite, then dim$(X)=I(W)$.\end{thm}
In the first part of this section we will explain how to build a cube complex associated to a given wall space. Next we will show that the complex thus obtained is CAT(0) and demonstrate that a proper action of a group on the wall space induces a proper action of the same group on this complex.

The vertices $X^0$ of the cube complex are given by a subset of the sections for the natural map
\begin{eqnarray*}\pi:\cH & \to & \cH/*=W\\
h &\mapsto & \{h,h^c\}.\end{eqnarray*}
Each section should be thought of as defining a transverse orientation on the set of walls which picks out, for each wall, a preferred ``side'', or half-space. The sections we require all satisfy the following compatibility condition:
$$\sigma(\oh)\not\subseteq\sigma(\ok)^c\hbox{ for any }\oh,\ok\in W.$$
In other words, if two walls $w_1$ and $w_2$ are disjoint, then the corresponding half-spaces $\sigma(w_i)$  are not disjoint. Notice that the condition says nothing for walls which cross, but for disjoint walls the following picture never occurs:
\begin{center}
\setlength{\unitlength}{1cm}
\begin{picture}(4,3)
\put(1,0){\line(0,1){2}}
\put(1,2.5){\makebox(0,0){$\oh$}}
\put(1,1){\vector(-1,0){1}}
\put(-0.7,1.3){\makebox(0,0){$\sigma(\oh)$}}
\put(3,0){\line(0,1){2}}
\put(3,2.5){\makebox(0,0){$\ok$}}
\put(3,1){\vector(1,0){1}}
\put(4.7,1.3){\makebox(0,0){$\sigma(\ok)$}}
\end{picture}\end{center}
Those sections satisfying this condition will be called \emph{admissible} and our vertex set is the set of all admissible sections. We connect two vertices with an edge if and only if, regarded as functions, their values differ on a single wall, so that now we have a graph. Each edge is labeled by the wall on which the associated sections differ. This graph is, in general, not connected, and we now explain which connected component to consider.

\smallskip 

Choose a point $p\in Y$ and define a section $\sigma_p$ for $\pi$ as follows: 
\begin{center}$\sigma_p(\oh)$ is the element of the set $\{h,h^c\}$ containing $p$.\end{center}
This is clearly a vertex of $X$, that we will refer to as \emph{special vertex}. We denote by $\Gamma_p$ the component of the graph containing $\sigma_p$.
\begin{lemma}For any $p,q\in Y$, then $\Gamma_p=\Gamma_q$ i.e., the vertices $\sigma_p, \sigma_q$ are connected by an edge path in the graph $\Gamma$.\end{lemma}
\begin{proof}First note that the sections $\sigma_p$ and $\sigma_q$ agree on all but the finitely many walls separating $p$ and $q$. Let $\{h_1, h_1^c\}, \ldots, \{h_n, h_n^c\}$ denote the set of walls on which they differ, and labeled so that for each $i$, $p\in h_i, q\in h_i^c$. The half-spaces $h_i$ form a finite partially ordered set (by inclusion) and so among them we can find a minimal element, which, by relabeling if necessary, we may assume to be $h_1$. Now consider the section $\sigma$ which is obtained from the admissible section $\sigma_p$ by setting $\sigma(\{h,h^c\})=\sigma_p(\{h,h^c\})$ if $h\not=h_1$ and $\sigma(\{h_1,h-1^c\})=h_1^c$, so that $\sigma$ agrees with $\sigma_p$ in every coordinate but one, where it agrees with $\sigma_q$ instead.
We claim that $\sigma$ is an admissible section. If this were not the case then $\sigma\{h, h^c\}=h, \sigma\{k, k^c\}=k$ for some walls $\{h, h^c\}, \{k, k^c\}$ with $h\subseteq k^c$. Since $\sigma$ agrees with the admissable section $\sigma_p$ in every wall except $\{h_1, h_1^c\}$ this must be one of the walls $\{h, h^c\}, \{k, k^c\}$, so there are two possibilities:

Case 1, $\{h_1, h_1^c\}= \{h, h^c\}$. We have $h_1=\sigma_p\{h_1, h_1^c\}\not=\sigma\{h_1,h_1^c\}=h_1^c=h$ and $k=\sigma_p\{k, k^c\}=\sigma\{k,k^c\}$. Hence $h_1^c\subseteq k^c$ or $k\subseteq h_1$ which does not contain $q$. The half-space $h_1$ was assumed to be minimal among those  containing $p$ but not $q$ and so $k$ cannot contain $p$ either. It implies that $k^c=\sigma_p\{k,k^c\}=\sigma\{k,k^c\}=k$ which is a contradiction.

Case 2, $\{h_1, h_1^c\}= \{k, k^c\}$. Here we have $h_1=\sigma_p\{h_1, h_1^c\}\not=\sigma\{h_1,h_1^c\}=h_1^c=k$ and $h=\sigma\{h, h^c\}=\sigma_p\{h,h^c\}$. Hence $h\subseteq k^c=h_1$  which does not contain $q$. The half-space $h_1$ was assumed to be minimal among those  containing $p$ but not $q$ and so $h$ cannot contain $p$ either. This implies that $h^c=\sigma_p\{h,h^c\}=\sigma\{h,h^c\}=h$ which is a contradiction.

Hence $\sigma$ is admissible and since it differs from $\sigma_p$ in a single coordinate it represents a vertex adjacent to $\sigma_p$. It also differs from $\sigma_q$ in only $n-1$ coordinates (the walls $\{h_2,h_2^c\},\ldots, \{h_n,h_n^c\}$). By induction on $n$ we obtain a sequence of vertices in $\Gamma$, $v_0=\sigma_p, v_1,\ldots, v_n=\sigma_q$ each adjacent to the next. This gives the required path from $\sigma_p$ to $\sigma_q$.\end{proof}
%
\begin{rem}\label{realization}
Corresponding to any special vertex $\sigma_p$ in $\Gamma_p$ there is a nonempty subset $x_{\sigma}$ of $Y$, which is the intersection of $\{\sigma_p(\oh)|\oh\in W\}$. The set $x_{\sigma}$ consists of all the points in $Y$ which are not separated from $p$ by any walls. It is reduced to a single point exactly when the wall space $Y$ is such that any two distinct points are separated by at least one wall. At that stage, if we say that two points in $Y$ are \emph{wall-equivalent} when there are no wall separating them, we can identify the set of special vertices in $X$ with the set of wall-equivalence classes in $Y$. There is an edge between two equivalence classes if and only if there is a single wall separating them. The pseudo-metric described in Section~\ref{basics} induces a metric on this subset of the vertices of $X$, which is equal to the combinatorial distance in  the graph $\Gamma_p$. Each edge in the connected component of an admissible vertex is labeled by a wall, and two half-spaces $h$ and $k$ have non-empty intersection exactly when they are in the image of a common section $\sigma$, i.e. when there exists $\sigma$ such that $\sigma(\oh)=h$ and $\sigma(\ok)=k$, where $\oh=\{h,h^c\}$ and $\ok=\{k,k^c\}$ are the walls associated to the half-spaces $h$ and $k$ respectively. In this case we will say that $\sigma$ \emph{lies in} $h\cap k$.\end{rem} 
We will now show how to attach cubes to the graph $X^1=\Gamma_p$ to construct a cubical complex. We first start with an easy fact.
\begin{lemma}\label{LoopsAreEven}Any closed loop has even edge length.\end{lemma}
\begin{proof} As noted above passing along an edge in $X^1$ corresponds to changing the value of a section on a single wall. As we proceed along a closed loop the value on each wall must be changed an even number of times to get back to the original section $\sigma$, so we must pass along an even number of edges.\end{proof}
For $k\in{\bf N}$, we will define a \emph{$k$-corner} to be a vertex $\sigma$ in $X^1$ together with a family $e_1,\dots,e_k$ of edges incident to $\sigma$ such that the walls $\oh_1,\dots,\oh_k$ labeling the edges pairwise cross. For each $i$ let $\tau_i$ denote the vertex adjacent to $\sigma$ along the edge $e_i$. We have the following.
\begin{lemma} Any $k$-corner $\sigma$ is
 contained in a unique subcomplex $C\subseteq X^1$ isomorphic to
 the $1$-skeleton of a $k$-cube.
\end{lemma}
\begin{proof}We can represent a $k$-dimensional cube $C^k$ as follows:
$$C^k=\{v:\{\oh_1,\dots,\oh_k\}\to [0,1]\},$$
so to prove the lemma it is enough to simplicially embed the 1-skeleton of $C^k$ in $X^1$, with one corner of $C^k$ mapped to $\sigma$. To do this, we first map the vertex $(0,\dots,0)$ to $\sigma$, and the vertices $v_i(\oh_j)=\delta_{ij}$ to $\tau_i$ for all $i=1,\dots,k$ (where $\delta_{ij}$ is the usual Kronecker symbol). For a generic vertex $v\in C^k$, we first extend $v$ to $W$ as follows:
\begin{eqnarray*}\hat{v}:W & \to & [0,1]\\
\oh &\mapsto &{\left\{\begin{array}{cc}v(\oh_i)& \hbox{ if }\oh=\oh_i \hbox{ for some $i$}\\
0 & \hbox{ otherwise }\end{array}\right.}
\end{eqnarray*}
and use this extension to define a map from the vertices $v$ of the cube $C^k$ to the set of  sections $W\to\cH$ as follows:
$$\tau_v(\oh)=\left\{\begin{array}{cc}\sigma(\oh)& \hbox{ if }\hat{v}(\oh)=0\\
\sigma(\oh)^c &  \hbox{ if }\hat{v}(\oh)=1.\end{array}\right.$$
We need to prove that all the sections constructed in this way are admissible. Denote by $|v|$ the number of non-zero coefficients of $v$. If $|v|=1$, then $\tau_v$ is adjacent to $\sigma$ via one of the edges $e_i$ by construction, hence it is one of the $\tau_i$, which is admissible by construction. Let us assume by induction that $\tau_v$ is admissible for any $|v|\leq n$ and some $n<k$. For $|v|=n+1$, label the non-zero values of $v$ by $\oh_1,\dots,\oh_{n+1}$. We know by the induction hypothesis that the section
$$\tau(\oh)=\left\{\begin{array}{ccc}\sigma(\oh)& \hbox{ if } & \hat{v}(\oh)=0\\
\sigma(\oh)^c &  \hbox{ if } & \oh=\oh_j, j=1,\dots,n\\
\sigma(\oh)^c &  \hbox{ if } & \oh=\oh_{n+1}\end{array}\right.$$
is admissible, so it remains to check that the section
$$\tau_v(\oh)=\left\{\begin{array}{cc}\tau(\oh)& \hbox{ if }\oh\not=\oh_{n+1}\\
\sigma(\oh)^c &  \hbox{ if }\oh\not=\oh_{n+1}\end{array}\right.$$
is admissible, i.e. that $\tau(\oh_{n+1})^c$ is never contained in $\tau(\oh)^c$ for $\oh_{n+1}\not=\oh$, but this is clear since if $\oh=\oh_s$ for $s=1,\dots,k$ and $s\not=n+1$ then $\oh_s$ crosses $\oh_{n+1}$, and if $\oh\not=\oh_s$ for all $s=1,\dots,k$, since the section obtained from $\sigma$ by $v(\oh_{j})=\delta_{n+1j}$ is admissible, then $\sigma(\oh_{n+1})=\tau(\oh_{n+1})^c=\tau_v(\oh_{n+1})^c$ cannot be contained in $\tau_v(\oh)^c=\sigma(\oh)^c=\tau(\oh)^c$. The embedding is simplicial because two vertices of $C^k$ are adjacent exactly if their defining maps differ on a single $\oh_i$, which means that the sections obtained differ on exactly a single wall, hence are adjacent in $X^1$. Notice that the embedding is completely determined by the corner itself, so that each corner defines a unique cube.
\end{proof}
We can now define a cube complex by gluing a $k$-dimensional cube to each $k$-corner. The existence of a $k$-dimensional cube implies the existence of a $k$-corner and hence $I(W)$ is the dimension of the cube complex when $I(W)$ is finite.

It remains to show that the cube complex is CAT(0). To do so we will use Gromov's theorem (Theorem \ref{Gromov} of this note) in the spirit of the Cartan-Hadamard theorem, that a geodesic cell complex is CAT(0) if and only if it is locally CAT(0) and simply connected. 
\begin{lemma}The cube complex $X$ is locally CAT(0).\end{lemma}
\begin{proof} That the link is simplicial is straightforward: A $k$-simplex in the link corresponds to a $(k+1)$-corner, and the intersection of two corners is again a corner (of smaller dimension).
That the link of a vertex is a flag complex is also more or less by construction, as the vertices of the 1-skeleton of a $(k-1)$-dimensional simplex in the link of a vertex $\sigma$ are labeled by hyperplanes $\oh_1,\dots,\oh_k$ pairwise sharing a square, i.e. pairwise crossing. This means that $\sigma$ is a $k$-corner, on which by construction has been glued a $k$-cube, and a $k$-cube gives a $(k-1)$-dimensional simplex in the link, from which we deduce that the 1-skeleton of a $(k-1)$-dimensional simplex actually belongs to a $(k-1)$-dimensional simplex.\end{proof}
\begin{lemma}The cube complex $X$ is simply connected.\end{lemma}
\begin{proof}Choose a vertex $\sigma_0$ in $X$ and let $\gamma$ be a loop in $X$ based at $\sigma_0$. Homotope $\gamma$ into the 1-skeleton $X^1$, and regard it as an edge path, which we will  denote by $\ell$. It suffices to show that this edge path can be homotoped to a point. Without loss of generality we can assume that $\ell$ does not backtrack. Let $\sigma$ denote any vertex on $\ell$  which maximizes the distance from $\sigma_0$, in the combinatorial distance $d_1$ on the $1$-skeleton of $X$. This is the cardinality of a shortest edge path between two points, which amounts to the cardinality of the set of walls on which the two sections differ, and which we will denote by $n$. 

Denote by $a$ and $b$ the two adjacent vertices to $\sigma$ in the loop $\ell$. Due to Lemma~\ref{LoopsAreEven}, combined with the assumption that $\sigma$ is furthermost from $\sigma_0$, we know that $d_1(\sigma_0,a)=d_1(\sigma_0,b)=n-1$. Denote by $\oh_a$ (resp. $\oh_b$) the wall on which the section defining $\sigma$ differs from the one defining $a$ (resp. $b$). Since $\ell$ does not backtrack we know that the vertices $a$ and $b$ are distinct so that the walls $\oh_a, \oh_b$ are distinct and carry distinct labels. We claim that $\sigma$ with the two edges to $a$ and $b$ is a 2-corner, i.e. that the two walls $\oh_a$ and $\oh_b$ cross, meaning that the four intersections of half-spaces $h_a\cap h_b$, $h_a^c\cap h_b$, $h_a\cap h_b^c$ and $h_a^c\cap h_b^c$ are all non-empty. Without loss of generality we can assume that $\sigma$ lies in $h_a\cap h_b$ (i.e. that $\sigma(\oh_a)=h_a$ and $\sigma(\oh_b)=h_b$, see end of Remark~\ref{realization}), and by assumption $a(\oh_a)=h_a^c$ and $b(\oh_b)=h_b^c$. But since $b$ only differs from $\sigma$ on the wall $\oh_b$ and $a$ only differs from $\sigma$ on $\oh_a$, and since $\oh_a\not=\oh_b$, $b$ agrees with $\sigma$ on the wall $\oh_a$, and hence $b(\oh_a)=h_a$ and similarly $a(\oh_b)=h_b$. Hence $b$ lies in $h_a\cap h_b^c$ and $a$ lies in $h_a^c\cap h_b$ so that they are both non-empty. It remains to show that $h_a^c\cap h_b^c$ is non-empty, which we do by showing that $\sigma_0$ lies in it. Suppose this were not the case. Then either $\sigma_0(\oh_a)=h_a$ or $\sigma(\oh_b)=h_b$ and by relabeling if necessary we may assume that $\sigma_0(\oh_a)=h_a$. It follows that the wall $\oh_a$ does not separate $\sigma$ and $\sigma_0$, so the set of walls on which $\sigma_0$ and $\sigma$ differ does not include $\oh_a$. But  the sections $\sigma$ and $a$ agree on all the walls on which $\sigma_0$ and $\sigma$ disagree, and also on the wall $\oh_a$ itself, hence $d_1(\sigma_0,a)=d_1(\sigma_0,\sigma)+1$ which is a contradiction and $\sigma_0$ lies in $h_a^c\cap h_b^c$ as required.

We have shown that there is a 2-corner at $\sigma$ so that  the following section is admissible:
$$\tau(\oh)=\left\{\begin{array}{ccc}\sigma(\oh)& \hbox{ if } & \oh\not=\oh_a,\oh_b\\
\sigma(\oh)^c &  \hbox{ if } & \oh=\oh_a\\
\sigma(\oh)^c &  \hbox{ if } & \oh=\oh_b\end{array}\right.$$
Furthermore $d_1(\sigma_0,\tau)=n-2$ (because $\sigma_0$ and $\tau$ coincide on exactly $n-2$ walls). We homotop $\ell$ to a closed loop $\ell'$  which coincides with $\ell$ everywhere except in a neighborhood of $\sigma$ where we replace the edges from $b$ to $\sigma$ and from $\sigma$ to $a$ by the edges from $b$ to $\tau$ and from $\tau$ to $a$. The homotopy is supported on  the unique square defined by these four vertices. Repeating this procedure to all furthermost vertices to $\sigma_0$ we obtain a loop homotopic to $\ell$ and lying in a ball of strictly smaller radius in the (integer valued) $d_1$ metric around $\sigma_0$. Iterating the procedure we will contract $\ell$ to the base point $\sigma_0$.
\end{proof}
If a group $G$ acts on the wall space $Y$, then it induces an action on the set of walls $W$, and we define an action on the vertices $X^0$ of the cube complex as follows:
$$g(\sigma)(\oh)=g\left(\sigma(g^{-1}(\oh))\right).$$
This action is obviously simplicial, because if $\sigma$ and $\tau$ share a common edge, it means that there is a unique wall $\ok$ on which they differ, so for any $g\in G$ and any wall $\oh$ such that $g^{-1}(\oh)\not=\ok$, then 
$$g(\sigma)(\oh)=g\left(\sigma(g^{-1}(\oh))\right)=g\left(\tau(g^{-1}(\oh))\right)=g(\tau)(\oh),$$ 
whereas for $g^{-1}(\oh)=\ok$, then  
$$g(\sigma)(\oh)=g\left(\sigma(g^{-1}(\oh))\right)=g\left(\sigma(\ok)\right)\not =g\left(\tau(\ok)\right)=g\left(\tau(g^{-1}(\oh))\right)=g\left(\tau(\oh)\right).$$ 
This means that $g(\sigma)$ and $g(\tau)$ differ on a single wall as well. Finally, if the action of $G$ on $Y$ is proper, then so is the action of $G$ on the cube complex $X$, precisely because the pseudo-metric given by the wall structure corresponds to the combinatorial distance given by the $1$-skeleton of $X$ (see Remark \ref{realization}). Notice that we didn't say anything about cocompactness, and indeed we shall see in the next section (Example~\ref{CocNonCoc}) that sometimes a cocompact action on a wall space yields a non-cocompact action on the corresponding cube complex.
\section{Examples}
\begin{ex}[Taken from \cite{N}] The Coxeter group $PGL_2({\bf Z})$ acts by isometries on the upper half-plane preserving a hyperbolic tessellation by isometric triangles, each with one ideal vertex (as illustrated in Figure~\ref{23inftycube} below). This gives the plane a wall space structure, where the walls are the mirrors for the reflections. The intersection number for this space with walls is 3, and the 3-cubes in the corresponding cube complex are given by the triple intersections of mirrors. The squares which are not faces of a 3-cube correspond to double points in the tessellation (i.e. intersections of two single walls). Part of the cube complex is also illustrated in Figure~\ref{23inftycube} and it may be viewed as a thickening of the Bass-Serre tree associated to the classical splitting of $\hbox{PGL}_2({\bf Z})$ as an amalgamated free product:
$$PGL_2({\bf Z})= D_2\mathop{*}_{{\bf Z}_2}D_3.$$
The dihedral groups $D_2$ and $D_3$ in the splitting are the two finite special subgroups $\langle  s_1,s_3\rangle$ and $\langle s_1,s_2\rangle$ respectively, in the standard Coxeter presentation 
$$PGL_2({\bf Z})=\langle s_1,s_2,s_3\mid s_1^2=s_2^2=s_3^2=(s_1s_3)^2=(s_1s_2)^3=1\rangle.$$
\begin{figure}
\resizebox{\textwidth}{!}{\includegraphics{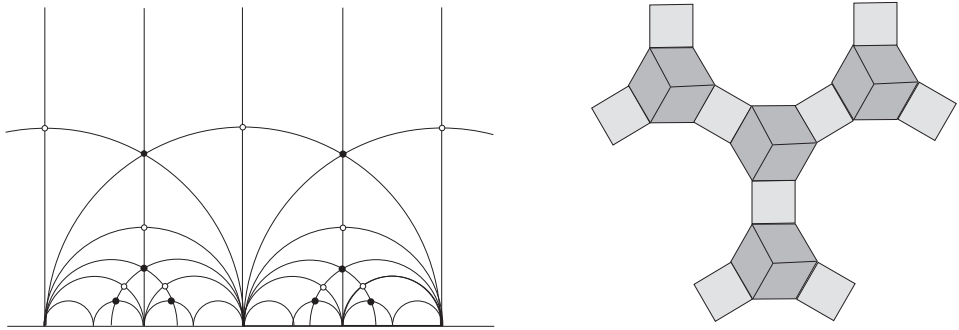}}
\caption{\label{23inftycube}The Coxeter complex for the group
$\text{PGL}_2({\bf Z})$ embedded in the hyperbolic plane, and the
corresponding CAT(0) cube complex.}
\end{figure}
\end{ex}
\begin{ex}\label{CocNonCoc}Consider the Euclidean Coxeter group 
$$G=\langle a,b,c\mid a^2=b^2=c^2=(ab)^3=(bc)^3=(ca)^3=1\rangle.$$ 
This group acts co-compactly on the Euclidean plane ${\bf E}^2$ as a reflection group, and the $G$-equivariant tiling by equilateral triangles defines a a discrete and $G$-equivariant family of geodesic lines which form the walls $W$. Let us explain how the wall space structure gives, under the construction described in Section~\ref{sageev}, a cubing which is isomorphic to that of ${\bf R}^3$ with its standard cubing. Pick $x_0$ a base point in ${\bf E}^2$, on the intersection of 3 lines that we label $O_x$, $O_y$ and $O_z$. Call ${\mathcal X}$, ${\mathcal Y}$ and ${\mathcal Z}$ the family of walls parallel to $O_x$, $O_y$ and $O_z$ respectively and choose $\omega$ a point at infinity which is not an end of $O_x$, $O_y$ or $O_z$. A vertex $\sigma$ of the cube complex is an admissible section, so it uniquely determines (due to compatibility conditions) a pair of consecutive parallel walls in each of the families ${\mathcal X}$, ${\mathcal Y}$ and ${\mathcal Z}$. For each of those 3 pairs of walls, choose the one which is nearest to $\omega$. The number of walls of ${\mathcal X}$, ${\mathcal Y}$ and ${\mathcal Z}$ respectively which are in-between those 3 walls and $O_x$, $O_y$ and $O_z$ (with a positive sign if the wall is in-between $O_x$, $O_y$ or $O_z$ and $\omega$ and a negative sign otherwise) gives an isometric isomorphism of the vertices of our cubing with ${\bf Z}^3$, which is the 0-skeleton of the standard cubing of ${\bf R}^3$.

The action of the group is affine, and preserves the plane perpendicular to the vector $(1,1,1)$ on which it acts co-compactly, however it cannot act co-compactly on ${\bf R}^3$ for dimension reasons. 

As pointed out by Valette, this group has Serre's property (FA), i.e., any action on a tree has a global fixed point. This follows from the observation that the group is generated by two elements of finite order whose product is also of finite order. Nonetheless our example shows that it acts cellularly and properly on a product of trees.\end{ex}
%

\end{document}